\documentclass[a4paper,english,final,titlepage,11pt]{article}
\usepackage[english]{babel}
\usepackage{amssymb}
\usepackage{amsmath}
\usepackage{indentfirst}
\usepackage{graphicx}
\usepackage{graphics}
\usepackage{epsfig}
\usepackage{multirow}
\usepackage{color}

\newcommand{\C}{\mathbb{C}}
\newcommand{\R}{\mathbb{R}}

\newcommand{\fimdem}{$\hfill \blacksquare$}
\newcommand{\e}{\mathrm{e}}
\def\diag{{\rm diag}}
\def\Tr{{\rm Tr}}

\newcommand{\beq}{\begin{equation}}
\newcommand{\eeq}{\end{equation}}

\newcommand{\bmat}{\begin{displaymath}}
\newcommand{\emat}{\end{displaymath}}

\newcommand{\vs}{\vspace{0.5cm}}

\def\1{{\bf 1}}

\newtheorem{theorem}{Theorem}[section]

\newtheorem{rmk}{Remark}[section]
\newtheorem{corol}{Corollary}[section]

\newtheorem{Exa}{Example}[section]

\begin{document}

\begin{center}
{\Large Matrices with  Hyperbolical Krein Space Numerical Range}

\vs

N. Bebiano\footnote{CMUC and  Mathematics Department,
 University of Coimbra,
3001-501 Coimbra (bebiano@mat.uc.pt)}, R. Lemos\footnote{CIDMA, Department of Mathematics, University of
Aveiro, Aveiro, Portugal (rute@ua.pt)} and G.
Soares\footnote{CM-UTAD, Department of  Mathematics, University of
Tr\'as-os-Montes  e Alto Douro,  
 Vila Real, Portugal (gsoares@utad.pt)}
 
 \end{center}


{\bf Abstract.} This paper is devoted to matrices with hyperbolical Krein space numerical range. This shape characterizes the 2-by-2 case and persists for certain classes of matrices, independently of their size. Ne\-cessary and sufficient conditions for low dimensional tridiagonal matrices to have this shape are obtained only involving the matrix entries.

\vs



{\bf Keywords;} numerical range, Krein space,  centrosymmetric matrices, tridiagonal matrices, Hyperbolical Range Theorem


\section{Introduction}

The investigation of the numerical range and its generalizations has attracted the interest of researchers from different areas of pure and applied science \cite{BLP2004,BP}, since the concept was born with Toeplitz \cite{Toep} and Hausdorff \cite{Haus}
in the second decade of the last century. The numerical range of linear operators in an indefinite inner product
space is a topic which has deserved the study by some authors; we refer the interested reader to \cite{BLPS, BNLPS2005, BLS_Hyp, bebiano, ricardo0, ricardo, LR PAMS, LR ELA3, CTU, PT}. Ilya Spitkovsky and collaborators published several inspiring papers on the ellipticity of the classical numerical range (see e.g.\! \cite{BPSK, BS, Geryba, JiangSpi, KeeRodSpi}). Recently, the numerical range of some structured tridiagonal  matrices was studied \cite{BLS_Kac,BLS_AOT}. In this note, we investigate classes of matrices with  hyperbolical Krein space numerical range, which can be seen as counterparts of these results for the indefinite inner product space context. 
 
 Let $M_n$ stand for the associative algebra of
$n \times n$ complex matrices and $I_n$ be the identity matrix.
Consider the complex vector space ${\bf C}^n$ equipped with the Krein space
structure  induced by the {\it indefinite inner product}  
$[x,y]_J:= y^*J\,x$, for $x, y \in {\bf C}^n$, with $J\!= P^T(I_r \oplus -I_{n-r})\,P$,  for some permutation matrix  $P\in M_n$ and $0\leq r\leq n$ (for references on Krein spaces see e.g.\! \cite{azizov, LR}).    
 
The {\it Krein space numerical range (KNR)},  also called {\it indefinite numeri\-cal range},  
of $A\in M_n$ is denoted and defined as
$$W^J(A):=-W^{J}_-(A)\,\cup\, W^J _{\,+}(A)$$
with
$$
  W_{\pm }^J(A):=\left\{[Ax,x]_J:\,x\in\C^n, \,[x,x]_J={\pm 1}\right\}.
$$

It is clear that $W^J_-(A)= W^{-J}_{\,+}(-A)$ 
 and $W^J(A)=\{1\}$ if and only if $A=I_n$. 
 If $J=I_n$, then $W^J(A)=W^J_+(A)$ reduces to the well-known classical {\it numerical
  range} $W(A)$,  a compact and convex set as asserted by the   Toeplitz-Hausdorff Theorem  \cite{Haus,Toep}. 
  
  These sets share some properties, but they  also have big differences,  such as the  possible degeneracy of $W^J(A)$
	to the complex plane, which never occurs for $W(A)$ when  $A\in M_n$.  
Although  each set $W^J_{\pm}(A)$ is connected, it may not be closed or bounded \cite{CTU}. Both sets $W_+^J(A)$, $W_-^J(A)$ are convex, and  $W^J(A)$  is {\it pseudo-con\-vex}, that is, for any pair of distinct  points $z_1,z_2\in W^J(A)$, either 
$W^J(A)$ contains the closed line segment 
$\{tz_1+(1-t)z_2: ~0\leq t\leq 1\}
$ 
or $W^J(A)$ contains the half-lines 
$
\{t z_1+ (1-t) z_2 : ~t\leq 0 ~ \hbox{or} ~ t\geq 1\}.
$

A matrix $A\in M_n$ is called $J$-{\it Hermitian} if it equals  its $J$-{\it adjoint} which is given by $A^\#:=JA^*J$
 and it is  $J$-{\it unitary} if $AA^\#=I_n$. 
Consider the following representation
$$
 A=\Re^J(A)+i\,\Im^J(A),
$$
where $$\Re^J(A) =\frac{A+A^\#}{2}\qquad \hbox{and}\qquad
\Im^J(A) =\frac{A-A^\#}{2i}$$
are $J$-Hermitian matrices.
For  each angle $\theta\in \R$, let 
\begin{equation}\label{Htheta}
H_\theta(A)=\,
\Re^J(A) \cos\theta+\Im^J(A)\sin\theta, \qquad \theta \in \mathbb{R},
\end{equation}
which is also  $J$-Hermitian, and so its eigenvalues are
real or occur in complex conjugate pairs  \cite{LR}. It is known that if  $H_\theta(A)$ has 
non-real eigenvalues, then $W^J(H_\theta(A))$ is the whole real line \cite[Proposition~2.1]{BNLPS2005}. To avoid
this trivial case, we assume throughout that all its eigenvalues are real. 
Let us define
$$
  \sigma^J_{\pm}(H_\theta(A)):=\{\lambda\in\R:\ \exists\, x\in\C^n,\ [x,x]_J=\pm 1, \ H_\theta(A)\,x=\lambda x\},
$$ 
where
\[
\sigma^J_+\big(H_\theta(A)\big)=\{\lambda_1(\theta),\ldots,\lambda_r(\theta)\}, \quad \lambda_1(\theta)\geq\cdots\geq \lambda_r(\theta),
\] 
\[ 
\sigma^J_-\big(H_\theta(A)\big)=\{\lambda_{r+1}(\theta),\ldots,\lambda_n(\theta)\}, \quad \lambda_{r+1}(\theta)\geq\cdots\geq \lambda_n(\theta).
\] 
We will be mainly concerned with the class ${\mathcal J}$ of matrices $H_\theta(A)$ with real eigen\-va\-lues, such that $\sigma_+^J(H_\theta(A))$ and $\sigma_-^J(H_\theta(A))$ do {\it not interlace}, that is, 
either 
$
 \lambda_r(\theta)>\lambda_{r+1}(\theta)  
$ 
or
$
 \lambda_n(\theta)>\lambda_1(\theta). 
$  
Notice that  if the eigenvalues of $H_\theta(A)$ interlace, then $W^J(H_\theta(A))$ is the whole real line too. 

We recall that a {\it support  line} of a convex set $S\subset \C$ is a line containing a boundary point of $S$ and defining two half-planes, such that one of them does not contain $S$. The support lines of $W^J(A)$ are by definition the support lines of the convex  sets  $-W_-^J(A)$ and $W_+^J(A)$.

Let $\Omega_A=(\eta,\xi)$ be an interval of angles $\theta$, with greatest possible dia\-meter, such that $H_\theta(A) \in {\mathcal J}$. 
For $\theta \in \Omega_A$, if $\lambda_r(\theta)>\lambda_{r+1}(\theta)$ holds, then the maximal eigenvalue in $\sigma^J_-(H_\theta(A))$ is denoted by $\lambda_{L} \big(H_\theta(A)\big)$  
and the minimal eigenvalue in $\sigma^J_+\big(H_\theta( A)\big))$ 
 by $\lambda_{R} \big(H_\theta(A)\big)$. 
It is known that $x=\lambda_{L} \big(H_\theta(A)\big)$ and $x=\lambda_{R} \big(H_\theta(A)\big)$ are supports  lines of $-W^J_-(A)$ and $W^J_+(A)$, respectively.  
Analogously,  if 
$\lambda_n(\theta)>\lambda_1(\theta)$ holds,  the maximal eigenvalue in $\sigma^J_+(H_\theta( A))$  is denoted by $\lambda_{L} \big(H_\theta(A)\big)$ and the  minimal eigenvalue in $\sigma^J_-\big(H_\theta(A)\big)$ by $\lambda_{R} \big(H_\theta(A)\big)$). Due to the convexity of $W^J_{\pm}(A)$, these families of eigenvalues, when  $\theta$ ranges over  $\Omega_A$, provide the characterization of these sets as described below. 

The {\it KNR generating polynomial} of $W^J(A)$ is defined as
$$p_{A}(z, \theta) \,=\, {\rm det}\,\big(\Re^J\!(A)\cos\theta\,+\,\Im^J\!(A)\sin\theta\,-\,z\,I_n\big)$$%
and plays a central role in the study of this set.
 Indeed, to a matrix $A\in M_n$, through the equation 
$$
  {\rm det}\,\big(u\,\Re^J\!(A)\,+\,v\,\Im^J\!(A)\,+\,w\,I_n\big)\,=\,0,
$$ 
 is associated an  algebraic curve of {\it class}  $n$ in homogeneous line
coordinates, called the {\it boundary generating curve} of
$W^J(A)$.  
 Its real part, denoted by $C^J(A)$, gene\-ra\-tes  $ W^J(A)$ as its {\it
pseudo-convex hull}, which is obtained as follows: for any two points $z_1=[Ax,x]_J$, $z_2=[Ay,y]_J$ in the boundary gene\-ra\-ting curve, take the line segment joining them if $[x,x]_J[y,y]_J>0$ and the two rays $\{tz_1+(1-t)z_2: t\leq 0$ or $t\geq 1\}$ if $[x,x]_J[y,y]_J<0$ (see e.g.\! \cite{BLPS} for more details). 

All the possible shapes of the numerical range $W^J(A)$ of  matrices  of size $3$ were classified in \cite{ricardo0,ricardo}  based on the factorability of the KNR genera\-ting polynomial $p_A.$ If this polynomial factors into three linear factors, then $C^J(A)$ consists of three points, which are the eigenvalues of $A$, and  
$W^J(A)$ is the pseudo-convex hull of the eigenvalues.
An\-oth\-er possibility is that  $p_A$ factors into a linear factor and a quadratic one. Thus $C^J(A)$
consists of a point, the eigenvalue  $p$ corresponding to the linear factor, and a  hyperbola  or ellipse ${\mathcal H}$,  corresponding to the quadratic factor, according to the signature of the involved submatrix of $J$.
 If $p$ lies outside the hyperbola (resp.\! ellipse), then  
$p$ is a nondifferentiable boundary point and there are flat portions on the boundary of $W^J(A)$. 
If $p$ lies in the {\it interior} of ${\mathcal H}$, then $W^J(A)$ is the hyperbolic (resp.\! elliptic) disc bounded by ${\mathcal H}$ 
and the matrix $A$ may, or may not, be $J$-unitarily reducible.
We recall that a matrix $A$ is {\it $J$-unitarily reducible} if there exists a $J$-unitary matrix $U$, such that 
$$U^{\#} A\,U\,=\,A_1\oplus A_2,$$
 where the sum of the sizes of the blocks $A_1$, $A_2$ is the size of $A$. 

In general, not too much is known about the boundary of $W^J(A),$ denoted by $\partial W^J(A).$  Here, we restrict our attention to structured matrices with lower dimensional order. The investigation of the hyperbolicity of the Krein space numerical range of tridiagonal  matrices of order $3$ and of centrosymmetric tridiagonal  matrices of order $5$  with biperiodic main diagonal is investigated and the obtained results are given in terms of the matrices entries.  
We also consider centrosymme\-tric  tridiagonal matrices of order $4$ with {\it quasi}-biperiodic main diagonal, and we show that $C^J(A)$ may be formed by  two hyperbolas. In this case, flat portions may eventually appear on $\partial W^J(A)$ in conjunction  with hyperbolic arcs, or $W^J(A)$ may degenerate into the whole complex plane.  

The remaining of this note is organized as follows. Section~\ref{Prereq} contains auxiliary results useful in subsequent discussions. In Section~3, 
counterparts of  classical numerical range results, some recently given in \cite{BLS_AOT} and  others  mainly  due to 
Spitkovsky and collaborators  (e.g. \cite{BPSK, BS, KeeRodSpi}), are presented in the set up of indefinite inner product spaces. 
 Illustrative examples of the obtained results are also given. In Section~4, some observations on the case of $6\times 6$ matrices are performed.

\section{Pre-requisites}\label{Prereq}

For $A \in M_n$ the following basic properties of $W^J(A)$ will be used throughout. 
\begin{itemize}
  \item[(i)] $W^J(\alpha A+\beta I_n)\!=\alpha W^J(A)+\beta$ for any $\alpha, \beta\in \C$.
  \item[(ii)]   $W^J(A)$ is $J$-{\it unitarily invariant}: 
      $W^J(U^\#AU)\!=\!W^J(A)$  for any $J$-uni\-tary $U\in M_n$ . 
   \item[(iii)]  $W^{J[\kappa ]}\big(A[\kappa ]\big) \subseteq W^J(A)$, where $A[\kappa ]$ (resp.\! $J[\kappa]$) denotes the principal submatrix of $A$ (resp.\!\! $J$) that lies in rows and columns indexed by the  set $\kappa\subset\{1, \dots, n\}$     \cite[Lemma 2.3]{GLS}.
\item[(iv)] For $A=(a_{ij})$  tridiagonal, $W^J(A)$ is invariant under interchange of the entries $a_{j,j+1}$, $a_{j+1,j}$, $j=1,\dots,n-1$ \cite[Lemma 2.1]{bebiano}.
			 \item[(v)]{\it Hyperbolical Range Theorem}  \cite{BLPS, BLS_Hyp}: For $J=\diag(1,-1)$ and  $A\in M_2$, the set  
  $W ^J(A)$ is bounded by a non-degenerate hyperbola with foci at the eigenvalues $\lambda_1,\lambda_2$ of $A$, transverse and non-transverse axes of length
$$
 \left(\Tr(A^\#A) -2\Re(\lambda_1\bar\lambda_2)\right)^{\frac{1}{2}}\quad \hbox{and}\quad \left(|\lambda_1|^2+|\lambda_2|^2-\Tr(A^\#A)\right)^{\frac{1}{2}},
$$
respectively, if and only if
$$ 
  2\Re(\bar\lambda_1\lambda_2)<\Tr(A^\#A)< |\lambda_1|^2+|\lambda_2|^2.
$$
If $A\in M_2$ is  $J$-Hermitian non-scalar, then $W^J(A)$ degene\-rates into either two half-rays with endpoints at the eigenvalues $\lambda_1,\lambda_2$ of $A$, when the eigenvalues are real,  or the whole real line, when they are complex; and
it  reduces to a singleton if and only if $A$ is a scalar matrix. Moreover, $W^J(A)$ may be the whole complex plane, eventually  except a line (see \cite{BLPS SIAM Proc} for a detailed description).

 \item[(v)]If $z$ is a {\it corner} of $W^J(A)$, that is, if $z$ lies on more than one support line of $W^J(A)$, then $z$ is an eigenvalue of $A$ \cite{LR PAMS}.
\end{itemize}

\medskip

The following criterion  of non-degenerate hyperbolicity of $W^J(A)$ was obtained in \cite{BLS_Hyp}.

\begin{theorem}\label{T2.11}
Let $\widetilde{a},\widetilde{b}>0$, $J=I_{r}\oplus -I_{n-r}$, $0<r<n$, and $A\in M_n$. 
The set $W^J(A)$ is bounded by  the non-degenerate hyperbola centered at the origin, with horizontal transverse  
and vertical non-transverse semi-axes 
of lengths, respectively, $\widetilde a$ and $\widetilde b$  if and only if
\[
 \lambda_{R}\big(H_\theta(A)\big)=\left(\widetilde a^2 - \widetilde c^2 \sin^2\theta\right)^{\frac{1}{2}}
  \   \hbox{and} \ \   
 \lambda_{L}\big(H_\theta(A)\big)=-\left(\widetilde a^2 - \widetilde c^2 \sin^2\theta\right)^{\frac{1}{2}}
\]
where $\widetilde c^2= \widetilde a^2+ \widetilde b^2$ for   all $\theta \in \Omega_A=(-\theta_0, \theta_0)$ and $\theta_0=\arctan(\widetilde a/\widetilde b)$. 
\end{theorem}

\smallskip

The criterion in Theorem 2.1 can be modified as follows. 

\smallskip

  \begin{corol}\label{corol Hyp} 
Let
$p,q,t \in \mathbb{R}$, 
   $s=q+it$,  $\gamma=\frac{1}{2}\arg(s)$, $p^2<q^2+t^2$ 
   and let ${\mathcal H}$ be the hyperbola centered at the origin,  transverse  semi-axis paralell to
$\e^{i\gamma}$ and  non-transverse semi--axis of lengths $\sqrt{|s|+p}$ and $\sqrt{|s|-p}$, respectively.
If $J=I_{r}\oplus -I_{n-r}$, 
$0<r<n$  and $A\in M_n$,
 $W^J(A)$ is bounded by the non-degenerate hyperbola ${\mathcal H}$ 
if and only if 
 $z= \lambda_{R}(H_\theta(A))$  and $z=\lambda_{L}(H_\theta(A))$
 satisfy the quadratic equation
\begin{equation}\label{quadr_eq0}
z^2\,=\, p+q\,\cos (2\theta)+t\,\sin (2\theta)
\end{equation}
for all \,$\theta \in \Omega_A\,=\,(\gamma -\theta_0, \gamma +\theta_0)$ and
$$\theta_0=\arctan  \left(\frac{|s|+p}{|s|-p}\right)^{\!\frac{1}{2}}.
$$ 
\end{corol}
\noindent{\bf Proof.}
Under the hypothesis, we have $|s|>|p|$. Remark that $\partial\,W^J(A)={\mathcal H}$ if and only if $\partial\,W^J({\rm e}^{-i\gamma}A)={\rm e}^{-i\gamma}\mathcal H$ is the  rotated hyperbola centered at the origin, with horizontal transverse   and vertical non-transverse semi-axes, res\-pec\-tively, of lenghts $\widetilde a=\sqrt{|s|+p}>0$ and $\widetilde b=\sqrt{|s|-p}>0$. By Theorem~\ref {T2.11},  this is equivalent to
\[
 \lambda_{R}\big(H_\eta({\rm e}^{-i\gamma}A)\big)=\left( p+|s|-2|s|\sin^2\eta\right)^{\frac{1}{2}}=-\lambda_{L}\big(H_\eta({\rm e}^{-i\gamma}A)\big)
\]
for  $\eta \in (-\theta_0, \theta_0)$ and $\theta_0=\arctan(\widetilde a/\widetilde b)$, observing that ${\widetilde a}^2+{\widetilde b}^2=2|s|$. 
Since $H_{\eta}({\rm e}^{-i\gamma}A)=H_{\eta+\gamma}(A)$, considering $\theta=\eta +\gamma$,
 the previous hyperbolicity condition holds if and only if 
$$\lambda_{R}(H_{\theta}(A))\,=\, \lambda_{R}(H_{\theta-\gamma}({\rm e}^{-i\gamma}A))\,= \big(p+|s|-2|s|\sin^2(\theta -\gamma)\big)^\frac{1}{2},$$ 
$$\lambda_{L}(H_{\theta}(A))\,=\,-\lambda_{R}(H_{\theta}(A)),\qquad \theta \in (\gamma-\theta_0,\gamma+\theta_0).$$
By trivial trigonometric transformations, we find
 \begin{eqnarray*}
|s| -2|s|\,\sin^2(\theta-\gamma) 
&\!=\!& |s|\,\cos^2(\theta-\gamma)-|s|\,\sin^2(\theta-\gamma)\\
&\!=\!& |s|\,\cos(2\theta-2\gamma)  \\
& \!=\!& |s|\,\cos(2\gamma)\cos (2\theta)+|s|\,\sin (2\gamma)\sin
(2\theta)  \\
&\!=\!&  q\,\cos (2\theta )+t\,\sin (2\theta).  
\end{eqnarray*}
Then the previous hiperbolicity condition occurs if and only if 
$z=\lambda_{R}(H_\theta(A))$ and $z=\lambda_{L}(H_\theta(A))$
satisfy condition (\ref{quadr_eq0}), having in mind that the RHS of the quadratic equation in  (\ref{quadr_eq0}) is positive for all angles $\theta$ in the interval $\Omega_A=(\gamma -\theta_0, \gamma +\theta_0)$.\fimdem\vs

\medskip
 In Theorem \ref{T2.11} and Corollary \ref{corol Hyp}, we may consider, for $0<r<n$, $J= P^T(I_{r}\oplus -I_{n-r})\,P$, with $P\in M_n$ a permutation matrix.
\smallskip

\begin{rmk} In parallel to the classical numerical range (see, e.g.\! \cite[Proposition 1]{JiangSpi}),
if $C^J(A)$ contains a hyperbola centered at a point $(c,d)$, then the KNR generating polynomial of $W^J(A)$ is divisible by a factor of the form
$$
(z-c \cos \theta {-} d \sin \theta)^2- \big(p+q\,\cos (2\theta )+t\,\sin (2\theta)\big), 
$$
where $p,q,t$ are defined as in Corollary \ref{corol Hyp}.
\end{rmk}

\section{Results}
 
\subsection{Hyperbolicity of $W^J(A)$ for orders $3$ and $5$}

\smallskip

In this subsection, we  will be concerned with the Krein space numerical range of 
 tridiagonal matrices  of odd order with biperiodic main diagonal, that is, matrices $A=(a_{ij})\in M_{2n+1}$, such that $a_{ij}=0$, $|i-j|>1$, with $a_{ii}=a_1$ if $i$ is odd and $a_{ii}=a_2$ if $i$ is even,   
 considering  
 $$J={\rm diag}(1,-1,1,\ldots,-1,1)\in  M_{2n+1}$$
with biperiodic main diagonal too. Without loss of generality,
  we may focus  on  tri\-dia\-gonal matrices with real biperiodic main diagonal, denoted by
  $T_{2n+1}({\bf c},{\bf a},{\bf b})$, where
$${\bf a}=(a,-a,a,\ldots, -a, a), \quad   {\bf b}=(b_1,b_2,\ldots,b_{2n}), \quad  
{\bf c}=(c_1,c_2,\ldots, c_{2n})$$ are 
the main diagonal,  the first upper and the first lower subdiagonals, res\-pectively. In fact, for $m=2n+1$, we easily see  that
\begin{equation}\label{A T}
A \, = \,{\rm e}^{i\tau}\, T_{m}({\bf c},{\bf a},{\bf b})+\delta\, I_{m},
\end{equation}  
considering   $a=\frac{1}{m}|\Tr(JA)|$, $\delta = \frac{1}{m}\Tr(A)$, $\tau =\arg\Tr (JA)$ and
\begin{equation}\label{bjcj}
b_j\,=\,a_{j, j+1}\,{\rm e}^{-i\tau},\qquad c_j\,=\,a_{j+1,j}\,{\rm e}^{-i\tau}, \qquad j=1, \dots, m-1.
\end{equation}
By property (i), $W^J(A)$ is readily obtained from $W^J\big(T_m({\bf c},{\bf a},{\bf b})\big)$. 
  
  \medskip\smallskip

We are interested in  deriving conditions ensuring  nondegenerate hyperbo\-li\-ci\-ty of the KNR. We start  with the $3 \times 3$ tridiagonal case. 
  \smallskip

\begin{theorem}  \label{lemma3}
Let $J=\diag(1,-1,1)$, 
$A=T_3 ({\bf c},{\bf a},{\bf b})$ with  ${\bf a}=(a,-a,a),$ ${\bf b}=(b_1, b_2)$,  ${\bf c}=(c_1,c_2)$ and $\Delta=a^2+b_1c_1+b_2c_2$.
The set
$W^J(A)$ is a nondegenerate hyperbolic disc  bounded by $\mathcal H$,
with foci at $\pm \Delta^{\frac{1}{2}}$ 
and non-transverse axis of length 
\begin{equation} \label{ntransv}
 \left(2|\Delta|-2a^2+|b_1|^2+|b_2|^2+|c_1|^2+|c_2|^2\right)^{\frac{1}{2}},
\end{equation}  
if and only if 
\begin{equation}\label{hypcond}
a^2-2|\Delta|< \Tr(A^\#A) < a^2+2|\Delta |.
\end{equation}
In this case, $C^J(A)$ is the union of the hyperbola $\mathcal H$ and the point $a$. 
\end{theorem}
 
 \noindent{\bf Proof.} Since $\det(A-\lambda I_3)=-(\lambda-a)(\lambda^2-\Delta)$, the eigenvalues of $A$ are easily obtained: $a,-\Delta^{\frac{1}{2}},\Delta^{\frac{1}{2}}$, and we 
find that 
\begin{eqnarray*}
H_{\theta} (A)\!&= & \frac{1}{2}\left[
\begin{array}{ccc}
 2a  & b_1-\overline{c_1} & 0  \\
 c_1-\overline{b_1} & -2a  &  c_2-\overline{b_2}  \\
 0 & b_2-\overline{c_2} & 2a  \\
\end{array}
\right]\cos \theta \\
 & + & \!\!\frac{1}{2i} \left[
\begin{array}{ccc}
 0  &  b_1+\overline{c_1} &0  \\
 c_1+\overline{b_1} & 0  &  c_2+\overline{b_2}  \\
 0 & b_2+\overline{c_2} & 0  \\
\end{array}
\right] \sin \theta.
\end{eqnarray*}
The KNR generating polynomial of  $W^{J}(A)$ is 
$$p_{A}(z,\theta)  = -(z-a \cos\theta) \left(z^2-p-q\cos (2\theta )-t\sin (2\theta )\right)$$
with
\begin{eqnarray*}
p & \!=\! & \frac{a^2}{2} -\frac{1}{4}\left(|b_1|^2+|b_2|^2+|c_1|^2+|c_2|^2\right) \, = \, \frac{1}{4}\big(\Tr(A^\#\!A) -a^2\big) ,\\
q & \!=\! & \frac{a^2}{2} +\frac{1}{2}\, \Re\, (b_1c_1+b_2c_2 ) \, = \, \frac{1}{2}\, \Re (\Delta),\\
t & \!=\! & \frac{1}{2}\,\Im\, (b_1c_1+b_2c_2 )  \, = \, \frac{1}{2}\, \Im (\Delta).
\end{eqnarray*}
Since $\cos (2 \theta)= 2\cos^2 \theta -1$, we can see that 
$$ q\cos (2\theta )+t\sin (2\theta ) \, =\,
 a^2\cos^2 \theta -\frac{a^2}{2} +\frac{1}{2}\,\Re \big({\rm e}^{-2i\theta}(b_1c_1+b_2c_2 )\big)\\
$$
and using $\Re \big({\rm e}^{-2i\theta}(b_1c_1+b_2c_2 )\leq |b_1c_1+b_2c_2|$, we get
\begin{equation}\label{pqr<acos}
 p+ q\cos (2\theta )+t\sin (2\theta ) \, \leq \, a^2\cos^2 \theta.
\end{equation}
We have $q+i\,t=\frac{1}{2}\Delta$ and  $p^2< q^2+t^2$ holds if and only if 
$$ 
\frac{1}{2}\big|\Tr(A^\#\!A) -a^2\big|< |\Delta|,
$$
which is equivalent to (\ref{hypcond}) as shown below.

Firstly, suppose that (\ref{hypcond}) holds.  
In this case, $\Delta\neq 0$, the LHS of (\ref{pqr<acos}) is positive  
and
 $p_{A}(z, \theta)=0$ if and only if  $z=a\cos  \theta $ or 
\begin{equation}\label{zH_thetaG}
z\, =\, \pm \big(
 p+ q\cos (2\theta )+t\sin (2\theta )
\big)^\frac{1}{2}
\end{equation}
for  $\theta\in (\Gamma-\theta_0,\Gamma+\theta_0)$, where $\Gamma=\frac{1}{2}\arg(\Delta)$ and 
 $$
  \theta_0=\arctan\left(\frac{|\Delta|+2p}{|\Delta|-2p}\right)^\frac{1}{2}.
 $$ 
Let $\lambda_\pm (\theta)$ be the eigenvalues of $H_{\,\theta}(A)$ given by  the  RHS of (\ref{zH_thetaG}).
We  find 
$
w_\theta=\left(b_2-\overline{c_2}\,{\rm e}^{2 i  \theta},0,\overline{b_1}\,{\rm e}^{2 i  \theta  }-c_1\right)
$
as an eigenvector of  $H_{\,\theta}(A)$  
 associated to the eigenvalue $a\cos \theta $, such that $[w_\theta,w_\theta]_{J}
 >0$.
Further,
$$ 
   u_{\theta,\pm}\,=\,\left(\overline{c_1}\,{\rm e}^{i \theta }-b_1 {\rm e}^{-i \theta}, 2(a \cos \theta - \lambda_{\pm}(\theta)), \overline{b_2}\,{\rm e}^{i \theta}-c_2 {\rm e}^{-i \theta}  \right) 
$$ 
is an eigenvector of  $H_{\,\theta}(A)$   associated with
$\lambda_\pm(\theta)$, 
satisfying
 \[
  [u_{\theta,\pm}, u_{\theta,\pm}]_{J} 
  \,=\,    8\lambda_\pm(\theta)\big(a\cos \theta-\lambda_\pm(\theta)\big)
\] 
 and so
 $[u_{\theta,+} u_{\theta,+}]_{J}
 \,[u_{\theta,-},u_{\theta,-}]_{J}
   < 0$, 
because 
$ H_{\,\theta}(A)$ is not diagonal, then  (\ref{pqr<acos}) does not hold as equality
and
\begin{equation}\label{neg}
  \big(\lambda_+(\theta)  -a  \cos \theta\big)\big(\lambda_-(\theta) -a  \cos \theta\big) > 0.
\end{equation}
Therefore,   one of the eigenvalues $\lambda_\pm(\theta)$ belongs to $\sigma^{J}_-\big( {H}_{\,\theta}(A)\big)$, the other  and $a\cos \theta$ belong to $\sigma^{J}_+\big(  {H}_{\,\theta}(A)\big)$. We conclude that 
${H}_{\,\theta}(A)\in {\mathcal J}$, with
 $$\lambda_{R}\big({H}_{\theta}({A})\big) \,=\,\lambda_+(\theta)
 \qquad \hbox{and}  \qquad \lambda_{L}\big({H}_{\theta}(A)\big)\,=\,\lambda_-(\theta) $$ 
 for all $\theta \in \Omega_{A}=(\Gamma-\theta_0, \Gamma+\theta_0)$.  

By Corollary \ref{corol Hyp},  $W^{J}(A)$ is a non-degenerate hyperbolic disc centered at the origin, with  transverse   axis parallel to ${\rm e}^{i\Gamma}$, that is, on the line containing the eigenvalues $\pm\Delta^\frac{1}{2}$, of length $(2|\Delta|+4p)^\frac{1}{2}$
  and with non-transverse  axis of length
$  (2|\Delta|-4p)^\frac{1}{2}$
 equal to (\ref{ntransv}). The foci of the hyperbola are $\pm \Delta^\frac{1}{2}$
and the remaining eigenvalue  $a$  of $A$ is  in the interior, or on the boundary, of this hyperbolic disc.

Conversely, if $W^{J}(A)$ is a non-degenerate hyperbolic disc with the  
foci and non-transverse axis given, we easily see that 
$$ 
\left(2|\Delta|-a^2 +\Tr(A^\#A) \right)^{\frac{1}{2}} \qquad \hbox{and} \qquad \left(2|\Delta|+a^2 -\Tr(A^\#A) \right)^{\frac{1}{2}}
   $$ 
are the axes lengths, this implying  that  (\ref{hypcond}) must be satisfied.\fimdem\vs

\begin{rmk}
The hyperbola describing   $\partial\,W^J(A)$ for $A$ in Theorem~\ref{lemma3}  under the  condition {\rm (\ref{hypcond})}  
can be characterized using the matrix invariants, namely, the eigenvalues $\lambda_1, \lambda_2, \lambda_3$ of $A$ and the trace of $A^\#A$. In fact, the non-transverse axis has length 
 $$
\left(|\lambda_1|^2+|\lambda_2|^2+|\lambda_3|^2-\Tr(A^\#A)\right)^\frac{1}{2},
$$
at least one of the eigenvalues equals the trace of $A$ and the remaining two are the foci of the hyperbola.
\end{rmk}

\begin{Exa}\label{ex1} 
Let $J=\diag(1,-1,1)$ and $A=T_3({\bf c},{\bf a},{\bf d})$ with  ${\bf a}\,=\,(4, -4, 4)$,  ${\bf b}\,=(3+i,1)$, ${\bf c}=(-4+i,-i)$. The eigenvalues of $A$, given by $4,-\sqrt{3-2 i}$ and $\sqrt{3-2 i}$,  are represented by green dots in Figure  1.  
The boundary generating curve $C^J(A)$ consists of the point $(4,0)$ and the hyperbola of Figure 1, being 
 $W^J(A)$ bounded by this hyperbola.

\begin{figure}[!h]
\centering
\includegraphics[scale=1.1]{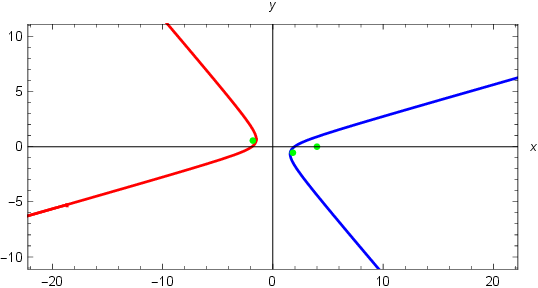}
\caption{ $C^J(A)$ for the matrices in Example \ref{ex1}}
\end{figure}
\end{Exa}

\newpage

Now, we apply the previous result to study the  case of centro\-symme\-tric tridiagonal matrices of order $5$.

\begin{theorem}\label{ordem5} Consider $J=\diag(1,-1,1,-1,1)$ and  $A=T_5({\bf c},{\bf a},{\bf b})$  with ${\bf a}=(a,-a,a,-a,a),$ ${\bf b}=(b_1, b_2,b_3,b_4),$ ${\bf c}=(b_4,b_3,b_2,b_1)$. Let
$$
\Delta_1=a^2+b_1b_4,\qquad M_1=2a^2-\left| b_1\right|^2-\left| b_4\right|^2,
$$
$$
\Delta_2=\Delta_1+2b_2b_3,\qquad M_2=M_1-2|b_2|^2-2|b_3|^2
$$
and ${\mathcal H}_k$ be the hyperbola with foci  $\pm\sqrt{\Delta_k}$ and non-transverse axis of length 
$$(2|\Delta_k|-M_k)^{\frac{1}{2}}.$$
The set
$W_J(A)$ is bounded by the nondegenerate hyperbola  ${\mathcal H}_2$ 
if and only if 
$|M_2| <2|\Delta_2|$.
In this case, $C_J(A)$ consists of  the point $a$ and the hyperbolas ${\mathcal H}_1$, ${\mathcal H}_2$, with ${\mathcal H}_1$ eventually degenerating into the points $\pm\sqrt{\Delta_1}$. 
\end{theorem}

\noindent{\bf Proof.} Let $E_2$ denote the exchange matrix of order $2$,  obtained from $I_2$ by reversing the order of its columns. Then
$$
Q=\frac{\sqrt{2}}{2}\left[\begin{matrix}
I_2 & 0 & E_2\\
0 & \sqrt{2} & 0\\
-E_2 & 0 & I_2
\end{matrix}\right]
$$
is such that $Q^\# Q=I_5$. Via this $J$-orthogonal matrix $Q$, the matrix $A$  is $J$-orthogonaly similar to the block diagonal matrix
$$
Q^{\#} A\, Q=\left[
\begin{array}{ccccc}
 a & b_1 & 0 & 0 & 0 \\
 b_4 & -a & 0 & 0 & 0 \\
 0 & 0 & a & \sqrt{2}\, b_3 & 0 \\
 0 & 0 & \sqrt{2}\, b_2 & -a & b_4 \\
 0 & 0 & 0 & b_1 & a \\
\end{array}
\right].
$$
Let $J_1=\diag(1,-1)$ and $J_2=\diag(1,-1,1)$. The boundary gene\-ra\-ting curves of $W^J(A)$ are those of $W^{J_1}(R)$ and $W^{J_2}(S)$, where  
$$
R=\left[
\begin{array}{cc}
 a & b_1 \\
 b_4 & -a \\
\end{array}
\right]\quad \hbox{ and } \quad S=\left[
\begin{array}{ccc}
 a & \sqrt{2} b_3 & 0 \\
 \sqrt{2} b_2 & -a & b_4 \\
 0 & b_1 & a \\
\end{array}
\right].
$$

Using properties (ii)-(iv), considering $\nu=\{2,3\}$, we easily find   
$$W^{J_1}(R) = W^{E_2J_1E_2}(E_2RE_2) = W^{J_2[\nu]}(S[\nu])  \subseteq W^{J_2}(S),$$ 
because $E_2RE_2=S[\nu]$.
The eigenvalues $\lambda_1$, $\lambda_2$ of $R$ are $\pm \sqrt{\Delta_1}$ and $$|\lambda_1|^2+|\lambda_2|^2 -\Tr(J_1R^*J_1R) =2|\Delta_1|-2a^2+| b_1|^2+| b_4|^2.$$ 
By the Hyperbolical Range Theorem, 
$W^{J_1}(R)$ is the hyperbolic disc bounded by ${\mathcal H}_1$. 

Now, applying  Theorem \ref{lemma3} to the tridiagonal matrix $S$ with the eigenvalues  $a,\pm\sqrt{\Delta_2}$,
we conclude that $W^{J_2}(S)$ is a nonde\-ge\-ne\-rate hyperbolic disc bounded by ${\mathcal H}_2$ 
 if and only if 
$$
a^2-2|\Delta_2| < \, \Tr(J_2\, S^* J_2 S) \, <a^2+2|\Delta_2|.
$$
Note that $$
\Tr(J_2\, S^* J_2 S)\,=\, 3 a^2-\left| b_1\right| {}^2-2 \left| b_2\right| {}^2-2 \left| b_3\right| {}^2-\left| b_4\right|^2,
$$
and the previous double inequality  may be witten as
$
 | M_2| <2|\Delta_2|.
$
 In this case, $a \neq 0$ and $C^{J_1}(R)={\mathcal H}_1$ may eventually degenerate into the pair of  points $\pm \sqrt{\Delta_1}$.
Then  $C^J(A)={\mathcal H}_1 \cup\, C^{J_2}(S)$ contains  also the point $a$ and the non-degenerate hyperbola ${\mathcal H}_2$.\fimdem\vs 

\begin{rmk} Under the hypothesis of Theorem \ref{ordem5}, we have
$$
\Tr(JA^* J\,A)\,=\,5 a^2-2 \left| b_1\right|^2-2 \left| b_2\right|^2-2 \left| b_3\right|^2-2 \left| b_4\right|^2.
$$
Considering   the eigenvalues $\lambda_k$ of $A$, $k=1, \dots, 5$, we easily find that
$$
\sum_{k=1}^5|\lambda_k|^2 -\Tr(JA^* J\,A) \,=\, \sum_{k=1}^2\left(2|\Delta_k|-M_k\right)
$$ 
is the sum of the squares of the non-transverse axes of the hyperbolas ${\mathcal H}_1$ and ${\mathcal H}_2$ belonging to $C^J(A)$. \end{rmk}

\begin{Exa}\label{ex2} 
Consider $J=\diag(1,-1,1,-1,1)$ and $A=T_5({\bf c},{\bf a},{\bf d})$  with ${\bf a}\,=\,(6, -6,6,-6,6)$, ${\bf b}=(-3,5,i,2)$,  ${\bf c}=(2,i,5,-3)$. The eigenvalues of $A$ given by $6$, $\pm\sqrt{30}$ and $\pm\sqrt{30+10 i}$ are represented by rose dots in Figure~2.
The boundary generating curve $C^J(A)$ consists of the point $(6,0)$ and two nested hyperbolas,  displayed in red and blue in Figure 2, being $W^J(A)$ bounded by the outer hyperbola.
\begin{figure}[!h]
\centering
\includegraphics[scale=1]{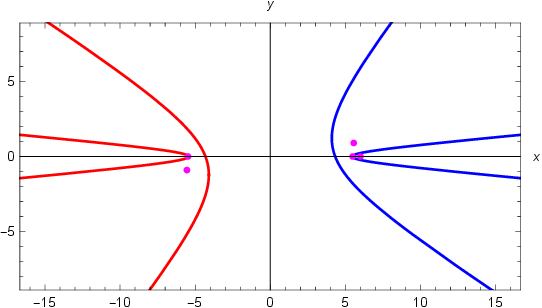}
\caption{$C^J(A)$ for the matrices in Example \ref{ex2}}
\end{figure}
\end{Exa}

\newpage

\subsection{Bihyperbolic $C^J(A)$ for matrices of order $4$}

\smallskip

In this subsection, we will consider 
 tridiagonal matrices of even order  $A=(a_{ij})\in M_{2n}$  
 with main diagonal, satisfying  $a_{ii}=a_{2n-i+1,2n-i+1}$ and $a_{ii}=a_1$ if $i\leq n$ is odd, 
 $a_{ii}=a_2$ if $i\leq n$ is even.  
 Let
\begin{equation}\label{Jeven}
J={\rm diag}(1,-1,\ldots,-1,1)\in  M_{2n},
\end{equation}
that is, with the main diagonal entries $j_1, \dots, j_n$ of $J$, satisfying $j_i=j_{n-i+1}$, $i=1, \dots, n$.
 Without loss of generality, by property (i), we may focus  on  tri\-dia\-gonal matrices with real {\it quasi}-biperiodic main diagonal, denoted by
  $T_{2n}({\bf c},{\bf a},{\bf b})$, where
$${\bf a}=(a,-a,\ldots, -a, a), \ \  {\bf b}=(b_1,b_2,\ldots,b_{2n-1}), \ \  
{\bf c}=(c_1,c_2,\ldots, c_{2n-1})$$ are 
the main diagonal, the first upper and  lower subdiagonals, respectively. In fact, let $m=2n$. Now, 
(\ref{A T}) holds too,  
considering   $a,\delta,\tau$ defined as before, $b_j, c_j$ written as (\ref{bjcj}), but with $J$ given by (\ref{Jeven}).

\smallskip

Now, the KNR of centrosymmetric tridiagonal matrices of order $4$ is studied. As shown below, 
$C^J(A)$  may be  bihyperbolic.  

\smallskip

\begin{theorem}\label{ordem4}  Let $J=\diag(1,-1,-1,1)$, $A=T_4({\bf c},{\bf a},{\bf b})$, ${\bf a}=(a,-a,-a,a)$,  ${\bf b}=(b_1,b_2,b_3)$ and   ${\bf c}=(b_3,b_2,b_1)$. Consider 
$$
\Delta =\, a^2 + b_1 b_3,
    \qquad 
\Delta_{\pm}=
     \Delta \mp \, a\, b_2 + \frac{b^2_2}{4},
$$
$$
M =\, 2 a^2-\left| b_1\right|^2-\left|b_3\right| ^2, \qquad M_\pm =\,M + \frac{1}{2}\left|b_2\right|^2\pm  2\,a\, \Re (b_2).
$$
\begin{itemize}
\item[\bf (a)] $C^J(A)$ is the union of two nondegene\-rate and nonconcentric  hyperbolas ${\mathcal H}_+$ and ${\mathcal H}_-$
if and only if $$
b_2 \neq 0, \qquad |M_+|<2|\Delta_+| \qquad \hbox{and}\qquad  |M_-| <2|\Delta_-|.
$$
In this case, $W^J(A)$ is the pseudo-convex hull of ${\mathcal H}_+$ and ${\mathcal H}_-$, where $\mathcal{H}_+$ $(\hbox{resp.} \ \mathcal{H}_-)$  
has  foci  
 $\frac{1}{2}\, b_2\pm\sqrt{\Delta_{+}}$ $\big($resp.\! $-\frac{1}{2}\, b_2\pm \sqrt{\Delta_{-}}\big)$ and 
 non-transverse axis of  length 
 $$
 \left(2|\Delta_+|-M_+\right)^{\frac{1}{2}} \qquad \left(\hbox{resp.} \ \big(2|\Delta_-|-M_-\big)^{\frac{1}{2}}\right).
 $$ 
 Flat portions on the boundary may  eventually exist, or $W_J(A)$ may degenerate into the whole complex plane.
\item[\bf (b)] 
 In particular, if 
$b_2=0$ and $|M| <2|\Delta|,$ then
 $W^J(A)$ is a nondegene\-rate hyperbolic disc 
with foci at $\pm\sqrt{\Delta}$,  non-transverse axis of  length 
 $
 \left(2|\Delta|-M\right)^{\frac{1}{2}}$.

\end{itemize}
\end{theorem}

\noindent{\bf Proof.} It is easy to see that  $A$ is $J$-orthogonally similar to the block diagonal matrix
$$
Q^{\#} A\,Q=\left[
\begin{array}{cccc}
 a & b_1 & 0 & 0 \\
 b_3 & -b_2-a & 0 & 0 \\
 0 & 0 & b_2-a & b_3 \\
 0 & 0 & b_1 & a \\
\end{array}
\right]$$
via the real $J$-orthogonal matrix  
$$
Q=\frac{\sqrt{2}}{2}\left[\begin{matrix}
I_2 & E_2\\
-E_2 & I_2
\end{matrix}\right]
$$ 
with $E_2$  the exchange matrix of order $2$.
Then  $C^J(A)$ is the union of the boundary generating curves of  $W^{J_1}(S_-)$ and $W^{-J_1}(R)$, with $J_1 = \diag(1,-1)$ and
$$
S_{\pm}  = \left[
\begin{array}{cc}
 a & b_1 \\
 b_3 & \pm b_2-a \\
\end{array}
\right],\qquad
R= \left[
\begin{array}{cc}
 b_2-a & b_3 \\
 b_1 & a \\
\end{array}
\right].
$$
Concerning the last set, it is clear that  $$W^{-J_1}(R)=W^{E_2J_1E_2}(E_2S_+E_2)=W^{J_1}(S_+).$$
By the Hyperbolical Range Theorem, each set
$W^{J_1}(S_\pm)$ is bounded by a non-degenerate hyperbola ${\mathcal H}_\pm$ with foci at the eigenvalues  of $S_\pm$,  
$$\lambda_{1,\pm}=\,{\pm}\frac{1}{2}\, b_2+\sqrt{\Delta_\pm}\,,\qquad  \lambda_{2, \pm}=\,\pm\frac{1}{2}\, b_2- \sqrt{\Delta_\pm}\,,$$ 
and non-transverse axes  of length given by the square root of
$$|\lambda_{1, \pm}|^2+|\lambda_{2, \pm}|^2 -\Tr(J_1\, S_\pm^{*}J_1\,S_\pm)
=2|\Delta_\pm|-M_\pm
$$ 
if and only if 
$
 -2|\Delta_\pm|< M_\pm  < 2|\Delta_\pm|.
$
The hyperbolas ${\mathcal H}_+$ and ${\mathcal H}_-$  are not concentric, except when $b_2=0$, in which case they coincide and (b) is obtained. Otherwise, $W^J(A)$ has a nondegenerate bihyperbolic boundary generating curve as described in (a). Henceforth, we may obtain:
\begin{itemize}
\item[i.] two nested hyperbolas as $C^J(A)$ and $W^J(A)$ is a hyperbolic disc bounded by the outer hyperbola; 
\item[ii.] hyperbolic arcs from ${\mathcal H}_+$  and ${\mathcal H}_-$ in $\partial W_J(A)$,  when these hyperbolas are not nested, and eventually some flat portions on  $\partial W^J(A)$, generated when taking the pseudo-convex hull of ${\mathcal H}_+$ and ${\mathcal H}_-$;
\item[iii.] the set $W^J(A)$ as the whole complex plane, when the eigenvalues of $H_\theta(A)$ interlace for some angles $\theta\in \Omega_A$. 
\end{itemize}
All these three cases may occur as illustrated in Examples \ref{ex3} - \ref{ex5} presented below.\fimdem

As the next example shows, the converse of Theorem \ref{ordem4} (b) does not hold, since $W^J(A)$ may be a hyperbolic disc even when $b_2\neq 0$.

\begin{Exa}\label{ex3} 
Let $J=\diag(1,-1,-1,1)$, $A=T_4({\bf c},{\bf a},{\bf d})$, ${\bf a}\,=\,(5, -5, -5,5)$,  ${\bf b}\,=(1,2,7)$ and ${\bf c}=(7,2,1)$. The eigenvalues of $A$, which are given by  $-1\pm\sqrt{43}$, $1\pm\sqrt{23}$,
are represented by green dots in Figu\-re~3. We have $C^J(A)=C^J(S_+)\cup C^J(S_-)$, where
 $$
S_{+}  = \left[
\begin{array}{cc}
 5 & 1 \\
 7 & -3 \\
\end{array}
\right]\qquad \hbox{and}\qquad S_{-}  = \left[
\begin{array}{cc}
 5 & 1 \\
 7 & -7 \\
\end{array}
\right],$$
using the notation in the proof of Theorem \ref{ordem4}.
Then $C^J(A)$ consists of  two nonconcentric nested hyperbolas displayed in Figure~3, where the blue (red) curves represent the branches of the hyperbolas that bound  $W^{J_1}_+(S_\pm)$ (resp., $-W^{J_1}_-(S_\pm)$). 
We conclude that $W^J(A)=W^J(S_+)$ is bounded by the outer hyperbola, centered at the point $(1,0)$ with focal distance equal to $2\sqrt{23}$ and vertical non-transverse axis of length $\sqrt{58}$. This example  illustrates the case i.\! presented in the proof above. 

\begin{figure}[!h]
\centering
\includegraphics[scale=0.8]{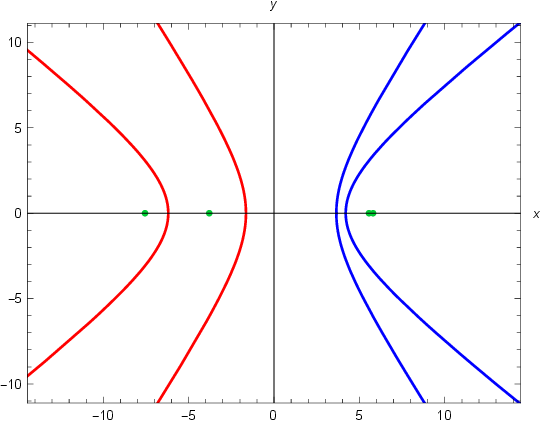}
\caption{ $C^J(A)$ for the matrices in Example \ref{ex3}}
\end{figure}
\end{Exa}

\newpage
\begin{Exa}\label{ex4} 
Let $J=\diag(1,-1,-1,1)$, $A=T_4({\bf c},{\bf a},{\bf d})$,  ${\bf a}\,=\,(6, -6, -6,6)$, ${\bf b}=(i,2+i,2)$ and  ${\bf c}=(2,2+i,i)$. The eigenvalues of $A$ are
$$
\frac{1}{2} \left((-2-i)\pm\sqrt{195+36 i}\right),\qquad \frac{1}{2} \left((2+i)\pm\sqrt{99-12 i}\right),
$$  
 represented by green dots and 
$C_J(A)$ consists of a pair of displaced hyperbolas displayed in Figure 4.
As in the previous example, the blue (red) curves represent the branches of the hyperbolas associated to points $z=[Ax,x]_J$ in  $C^J(A)$ generated by vectors $x$ with $[x,x]_J>0$ (resp., $[x,x]_J<0$).
In this case,  from the pseudo-convex hull of these two hyperbolas, flat portions will appear on  $\partial W^J(A)$. This example illustrates the case ii.\! in the proof of Theorem \ref{ordem4}.

\begin{figure}[!h]
\centering
\includegraphics[scale=1.3]{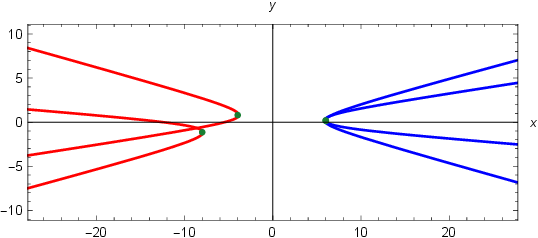}
\caption{ $C^J(A)$ for the matrices in Example \ref{ex4}}
\end{figure}

\end{Exa}

\begin{Exa}\label{ex5} 
Let  $J=\diag(1,-1,-1,1)$, $A=T_4({\bf c},{\bf a},{\bf d})$, ${\bf a}\,=\,(1, -1, -1,1)$, ${\bf b}=(3,5,1/2)$ and  ${\bf c}=(1/2,5,3)$. The eigenvalues of $A$ are
$\frac{1}{2} \left(-5\pm\sqrt{55}\right)$, 
$\frac{1}{2} \left(5\pm\sqrt{15}\right)$,
 represented by violet dots  and 
$C_J(A)$ consists of the pair of nonconcentric hyperbolas  displayed in Figure 5.
The blue (red) curves represent again the branches of the hyperbolas obtained from points $z=[Ax,x]_J$ in  $C^J(A)$ with generating vectors $x$ such that $[x,x]_J>0$ (resp., $[x,x]_J<0$).
In this case, we can easily conclude  that $W^J(A)$  is the whole complex plane.

\begin{figure}[!h]
\centering
\includegraphics[scale=1]{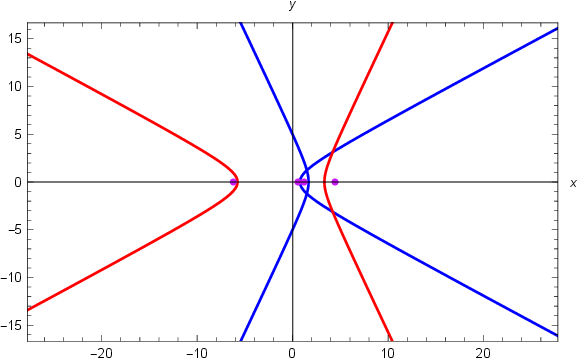}
\caption{$C^J(A)$ for the matrices in Example \ref{ex5}}
\end{figure}
\end{Exa}

\newpage

\section{Final remarks}

In the case of matrices of order $6$, the previous cenario changes.
Now, let $J=\diag(1,-1,1,1,-1,1)$  and   $A=T_6({\bf c},{\bf a},{\bf b})$ be  the centrosymmetric matrix with  ${\bf a}=(a,-a,a,a,-a,a),$ ${\bf b}=(b_1, b_2,b_3,b_4,b_5)$, ${\bf c}=(b_5,b_4,b_3,b_2,b_1)$. 
It may be seen that  $A$  is $J$-orthogonally similar to the block diagonal matrix: 
$$
Q^{\#} A\, Q\,=\left[
\begin{array}{cccccc}
 a & b_1 & 0 & 0 & 0 & 0 \\
 b_5 & -a & b_2 & 0 & 0 & 0 \\
 0 & b_4 & a-b_3 & 0 & 0 & 0 \\
 0 & 0 & 0 & a+b_3 & b_4 & 0 \\
 0 & 0 & 0 & b_2 & -a & b_5 \\
 0 & 0 & 0 & 0 & b_1 & a \\
\end{array}
\right],
$$
via the real $J$-orthogonal matrix   
$$
Q=\frac{\sqrt{2}}{2}\left[\begin{array}{cc}
I_3  & E_3\\
-E_3 & I_3
\end{array}\right].
$$
where $E_3$ is the exchange matrix of order $3$.

\smallskip

If $b_3=0,$ then  $C^J(A)=C^{J_2}(S)\,\cup\, C^{J_2}(R')$, where $J_2 = \diag(1,-1,1)$,
$$
S\, = \left[
\begin{array}{ccc}
 a & b_1 & 0  \\
 b_5 & -a & b_2  \\
 0 & b_4 & a
 \end{array}
\right]=\,E_2 R'E_2,\qquad
R' \, = \left[
\begin{array}{ccc}
a & b_4 & 0 \\
 b_2 & -a & b_5 \\
 0 & b_1 & a
\end{array}
\right].
$$
As $J_2=E_3J_2E_3$, $W^{J_2}(R')=W^{J_2}(E_2 R' E_2)=W^{J_2}(S)$ and $C^J(A)$ reduces to $C^{J_2}(S)$.
By Theorem \ref{lemma3}, the set
$W^{J_2}(S)$ is a nondegenerate hyperbolic disc, bounded by the hyperbola ${\mathcal K}$
with foci at the points $$\pm (a^2+b_1b_5+b_2b_4)^{\frac{1}{2}}$$ 
and non-transverse axis of length 
\[ 
\left(2|a^2+b_1b_5+b_2b_4|-2a^2+|b_1|^2+|b_2|^2+|b_4|^2+|b_5|^2\right)^{\frac{1}{2}},
\]  
if and only if 
\[
\left|\Tr(J_2\, S^* \! J_2 S)-a^2 \right| \,<\, 2|a^2+b_1b_5+b_2b_4|.
\]  
In this last case,  $C^{J}(A)=\{a\}\cup {\mathcal K}$  with the point $a$ belonging to the interior,  or the boundary, of the hyperbolic disc bounded by $\mathcal K$. 
 
 \smallskip

As it is well-known, tridiagonal matrices are too sensitive to perturbations, and even small perturbation provokes dramatic changes in their behaviour. In particular, the necessary and sufficient condition for hyperbolicity is violated for small changes in the  $(3,3)$  entry of $S$.

If $b_3\neq 0$, then  the boundary generating curves for both  tridiagonal blocks or order $3$ in $Q^{\#} A\, Q$ may not be hyperbolas. As illustrated in the next example, their KNR generating polynomials are irreducible cubics and the corresponding  boundary generating curves are algebraic curves of 6th degree, with hyperbolic-like shapes and deltoids in their interior.

\smallskip

\begin{Exa}\label{ex6} 
Consider $J=\diag(1,-1,1,1,-1,1)$ and $A=T_6({\bf c},{\bf a},{\bf b})$   with ${\bf a}=(4,-4,4,4,-4,4),$ ${\bf b}=(1, 2,-2,4,5)$, ${\bf c}=(5,4,-2,2,1)$. Moreover, let  $A_0=T_6({\bf c},{\bf a},{\bf b}_0)$, ${\bf b}_0=(1, 2,0,4,5)$, which only differs from the matrix $A$ on the entry $b_3$.  

From the previous discussion, it is clear that $W^J(A_0)$ is bounded by the hyperbola centered at the origin, with foci at the points 
$\pm\sqrt{29}$, horizontal tranverse axis and vertical non-transverse axis of lenght 
 $6\sqrt{2}$.  
This hyperbola is the dashed curve displayed in Figure 6, which forms the boundary generating curve of $W^J(A_0)$ together  with the point $(5,0)$  represented by a rose dot in its interior; the other two rose dots correspond to the remaining eigenvalues of $A_0$,  which are the foci given above.

\begin{figure}[!h]
\centering
\includegraphics[scale=0.9]{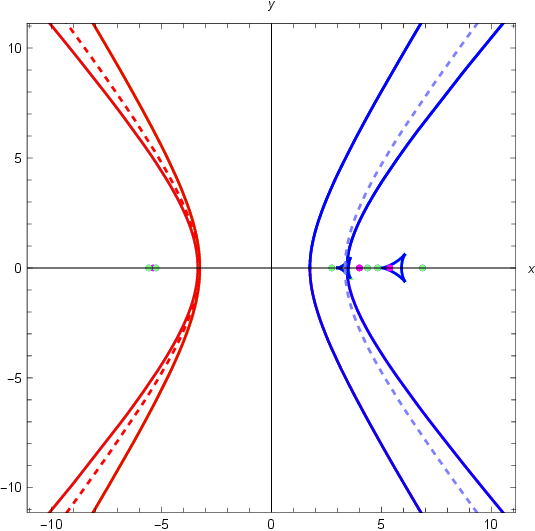}
\caption{$C^J(A_0)$ and $C^J(A)$ for the matrices in Example \ref{ex6}}
\end{figure}

Concerning the matrix $A$, the KNR generating polynomial of $W^J(A)$ is the product 
of two irreducible polynomials given by 
\begin{eqnarray*}
p_1(u,v,w)\! &=& \! w^3-68 u^3-11 u^2 w+3u\left(30 v^2+2 w^2\right)+18 v^2 w,\\
p_2(u,v,w)\! &=& \! w^3-20 u^3-11 u^2 w+2 u \left(27 v^2+w^2\right)+18 v^2 w.
\end{eqnarray*}
We can easily see that  $p_1(u,v,w)$ is an irreducible polynomial. In fact, supposing that
it is a product $p(u,v,w)q(u,v,w)$ of some real linear and quadratic homogeneous polynomials of type
\begin{eqnarray*}
p(u,v,w) &=& a u+b v+w\\
q(u,v,w) &=& g u^2+h u v+j u w+k v^2+l v w+w^2,
\end{eqnarray*}
after performimg this product and  using Mathematica 13.0, from the equality $p_1(u,v,w)=p(u,v,w)q(u,v,w)$, we obtain a contradiction.
Is is also clear that $p_1(u,v,w)$ is not the product of three real linear  homogeneous polynomials.
Analogousgly, we can prove that $p_2(u,v,w)$ is irreducible.
As displayed in Figure 6,  $C^J(A)$ contains two hyperbolic-like curves, which are not hyperbolas, and two deltoids in its interior. The left-most blue curve and right-most red curve represent the boundaries of the convex components $W^J_+(A)$ and $-W^J_-(A)$, respectively, that is, 
$W^J(A)$ is bounded that these two curves. Now, the eigenvaluess of $A$ are marked as the green dots.  

After several experimental tests,  we found that this type of shape  is kept for the Krein space numerical range of matrices of type $A'=T_6({\bf c},{\bf a},{\bf b'})$ with  ${\bf b'}=(1, 2,b_3,4,5)$, considering different values of $b_3$ in a neighbourhood of $0$,
and that the deltoids are getting smaller and smaller as we consider  values of $b_3$ approaching  $0$, indicating that we are aproaching hyperbolicty of the set, which only occurs for $b_3=0$.
\end{Exa}

\smallskip

\subsection*{Acknowledgment}

The  first author was partially supported by the Centre for Mathema\-tics of the University of Coimbra (funded by the Portuguese Government through FCT/MCTES, DOI 10.54499/UIDB/00324/2020). 
The se\-cond author was supported by the Center for Research and Development in Ma\-the\-ma\-tics and Applications (CIDMA) through the Portuguese Foundation for Science and Technology (FCT - Fundac\~ao para a Ci\^encia e a Tecnologia), project UIDB/04106/2020 (https://doi.org/10.54499/UIDB/04106/2020). 
The third author was financed by  CMAT-UTAD through the Portuguese Foundation for Science and Technology (FCT - Fundac\~ao para a Ci\^encia e a Tecnologia), reference UIDB/00013/2020, (https://doi.org/10.54499/UIDB/00013/2020).
 
 \smallskip


\end{document}